\documentclass[11pt]{article}
\usepackage{amsfonts, amssymb, amsmath}
\usepackage{enumerate}
\usepackage{esint}
\newtheorem{thm}{Theorem}[section]

\newtheorem{lem}[thm]{Lemma}
\newtheorem{prop}[thm]{Proposition}
\newtheorem{defn}[thm]{Definition}

\newcommand{\f}{\frac}

\newcommand{\vc}{\infty}

\newcommand{\RR}{\mathbb{R}}

\newcommand{\pa}{\partial}
\begin{document}
\title{New class of multiple weights and new weighted inequalities for multilinear operators}
\author{The Anh Bui\thanks{Department of Mathematics, Macquarie University, NSW 2109, Australia and Department of Mathematics, University of Pedagogy, HoChiMinh City, Vietnam. \newline{Email: the.bui@mq.edu.au and bt\_anh80@yahoo.com}  \newline
{\it {\rm 2010} Mathematics Subject Classification:} 42B20, 42B25, 35S05, 47G30.
\newline
{\it Key words:} Weighted norm inequality; Multilinear operator; Multilinear pseudodifferential; Commutator.}}

\date{}

\maketitle

\begin{abstract}
In this paper, we first introduce the new class of multiple weights $A^\vc_{\vec{p}}$ which is larger than the class of multiple weights in \cite{LOPTG}. Then, using this class of weights, we study the weighted norm inequalities for certain classes of multilinear operators and their commutators with new BMO functions introduced by \cite{BHS1}. Finally, we show that some multilinear pseudodifferential operators fall within the scope of the theory obtained in this paper.
\end{abstract}

 \tableofcontents
\section{Introduction and the main results}

The theory of multilinear Calder\'on-Zygmund singular integral operators, originated from
the work of Coifman and Meyer, has an important role in harmonic analysis. This direction of research has
been attracting a lot of attention in the last few decades, see for example \cite{CM1,CM2,CM3,GT, KS}
for the standard theory of multilinear Calder\'on-Zygmund singular integrals.

\medskip

Let $T$ be a multilinear operator initially defined on the $m$-fold
product of Schwartz spaces and taking values into the space of
tempered distributions,
\begin{equation*}
T: \mathcal{S}(\mathbb{R}^n)\times\ldots\times
\mathcal{S}(\mathbb{R}^n) \rightarrow \mathcal{S}'(\mathbb{R}^n).
\end{equation*}
By associated kernel to $T$ we shall mean the function $K$, defined off the diagonal $x = y_1 =\ldots=y_m$ in
$(\mathbb{R}^n)^{m+1}$, satisfying
\begin{equation*}
T (f_1, \cdots, f_m)(x) = \int_{(\mathbb{R}^n)^m} K(x,y_1,\ldots ,
y_m)f_1(y_1)\ldots f_m(y_m)dy_1 \ldots dy_m
\end{equation*}
for all $x \notin \cap^m_{j=1}$supp$f_j$.

\medskip

In this paper, we consider the following conditions:

(H1) For any $N>0$ there exists $C>0$ such that
\begin{equation}\label{CZ1}
|K(y_0, y_1, \ldots, y_m)|\leq
\frac{C}{(\sum_{k,l=0}^m|y_k-y_l|)^{N}}.
\end{equation}

(H2) For any $N>0$ there exists $C>0$ such that
\begin{equation}\label{CZ2}
|K(y_0,\ldots, y_j,\ldots,y_m)-K(y_0,\ldots, y'_j,\ldots,y_m)|\leq
\frac{C|y_j-y'_j|^\epsilon}{(\sum_{k,l=0}^m|y_k-y_l|)^{mn+\epsilon}}\min\{1, h^{-N}\},
\end{equation}
for some $\epsilon>0$ and all $0\leq j \leq m,$ whenever
$|y_j-y'_j|\leq \frac{1}{2}\max_{0\leq k \leq m}|y_j-y_k|:=h$.

(H3) There exist $1\leq q_1, \ldots, q_m<\vc$ and $1/q=1/q_1+\ldots+1/q_m$ such that $T$ maps continuously from $L^{q_1}\times\ldots\times L^{q_m}$ into $L^q$.

\medskip

It is clear that if $T$ satisfies (H1), (H2) and (H3) then $T$ falls within the scope of multilinear Calder\'on-Zygmund theory investigated by \cite{GT}. Therefore, according to \cite{GT}, if $1/p=1/p_1+\ldots+1/p_m$, the following statements hold:

(i) $T: L^{p_1}\times\ldots\times L^{p_m}  \rightarrow L^p$ when $1<p_1,\ldots,p_m<\vc$, and

(ii)  $T: L^{p_1}\times\ldots\times L^{p_m}  \rightarrow L^{p,\vc}$ when $1\leq p_1,\ldots,p_m<\vc$ and at least one $p_j=1$.

\medskip

The weighted norm inequalities of multilinear Calder\'on-Zygmund operators and their commutators with BMO functions were investigated in \cite{LOPTG}. In \cite{LOPTG}, the authors introduced the new maximal functions and multiple weights and then they proved that the new class of  multiple weights is suitable to study the boundedness of  multilinear Calder\'on-Zygmund operators and their commutators with BMO functions.

\medskip

Inspiring by the works of \cite{LOPTG} for boundedness of  multilinear Calder\'on-Zygmund operators and their commutators with BMO functions and of \cite{BHS1, BHS2} for the new class of weights and new BMO function spaces, the aim of this paper is to study the weighted norm inequalities of operators $T$ which satisfy (H1)-(H3) and their commutators by using the new BMO function spaces introduced by \cite{BHS1} and the new class of multiple weights introduced in Section 2.

\medskip

The organization of the paper is as follows. In Section 2, we introduce the new class of multiple weights and then investigate the weighted norm inequalities of some maximal functions. Section 3 establishes the main results of the paper. Firstly, the weighted estimates of multilinear operators $T$ are investigated (see Theorem \ref{thm1}). Secondly, we consider the weighted norm inequalities of the commutator $T_{\vec{b}}$ by using the new BMO functions and new class of multiple weights (see Theorem \ref{thm2}). In Section 4, we show that the obtained results can be applied to certain multilinear pseudodifferential operators.

\medskip

After finishing this paper, I had informed that the author in \cite{T} obtained the $A^\vc_p$ weighted norm inequalities for such an operator $T$. However,  in this paper, we study the the $A^\vc_{\vec{p}}$ weighted norm inequalities of $T$ and the obtained results and the new class of multiple weights in our paper are new. Moreover, the weighted norm inequalities of the commutator $T_{\vec{b}}$ with the new BMO functions in our paper are unique.

\section{Preliminaries}

To simplify notation, we will often just use $B$ for $B(x_B, r_B)$ and $|E|$ for the measure of $E$ for any measurable subset $E\subset \RR^n$.
Also given $\lambda > 0$, we will write $\lambda B$ for the
$\lambda$-dilated ball, which is the ball with the same center as
$B$ and with radius $r_{\lambda B} = \lambda r_B$. For each ball
$B\subset \RR^n$ we set
$$
S_0(B)=B \ \text{and} \ S_j(B) = 2^jB\backslash 2^{j-1}B \
\text{for} \ j\in \mathbb{N}.
$$

\subsection{The new class of weights and new BMO function spaces}
\subsubsection{Classes of multiple weights $A^\vc_{\vec{p}}$}
In this section, we would like to recall the definition of the new class of weights introduced by \cite{BHS2}. \\

For $1\leq p<\vc$ and $\theta\geq 0$, the weight $w$ ($w$ is a nonnegative and locally integrable function) is said to be in the class $A^{\theta}_{p}$ if there holds
\begin{equation}\label{classofnewweights}
\Big(\int_Bw\Big)^{1/p}\Big(\int_Bw^{-\f{1}{p-1}}\Big)^{1/p'}\leq C|B|(1+r_B)^{\theta}
\end{equation}
for all ball $B=B(x_B,r_B)$. In particular case when $p=1$, (\ref{classofnewweights}) is understood
$$
\f{1}{|B|}\int_Bw(y)dy\leq C(1+r_B)^{\theta}\inf_{x\in B}w(x).
$$
Then we denote $A^\vc_{p}=\cup_{\theta\geq 0}A^{\theta}_{p}$ and $A^\vc_\vc=\cup_{p\geq 1}A^\vc_p$.

We remak that $A^0_p$ coincides with the Muckenhoupt's class of weights $A_p$ for all $1\leq p<\vc$. However, in general, the class $A^\vc_p$ is strictly larger than the class $A_p$ for all $1\leq p<\vc$.
The following properties hold for the new classes $A^\vc_p$, see \cite[Proposition 5]{BHS2}.
\begin{prop}\label{property-Avc}
 The following statements hold:

i) $A_{p}^\vc\subset A_q^\vc$ for $1\leq p\leq q<\vc$.

ii) If $w\in A^\vc_p$ with $p> 1$ then there exists $\epsilon >0$ such that $w\in A_{p-\epsilon}^\vc$. Consequently, $A_p^\vc=\cup_{q<p}A_q^\vc$.

iii) If $w\in A^\vc_p$ with $p\geq 1$, then there exist positive numbers $\delta, \eta$ and $C$ so that for all balls $B$,
$$
\Big(\f{1}{|B|}\int_Bw^{1+\delta}(x)dx\Big)^{\f{1}{1+\delta}}\leq C\Big(\f{1}{|B|}\int_Bw(x)dx\Big)(1+r_B)^\eta.
$$
\end{prop}

In \cite{LOPTG}, to study the weighted norm inequalities of multilinear operators, the authors introduced the new maximal functions and the multiple weights. Adapting this idea to our situation, we introduce the new class of multiple weights.

In what follows, for given $m$ exponents $p_1, \ldots, p_m$, unless specified, otherwise we write $\vec{p}=(p_1,\ldots,p_m)$ and the number $p$ shall mean that
$$
\f{1}{p}=\f{1}{p_1}+\ldots+\f{1}{p_m}.
$$
For any number $r>0$, $r\vec{p}$ is defined by $r\vec{p}=(rp_1,\ldots,rp_m)$.

\begin{defn}
Let $1\leq p_1, \ldots, p_m<\vc$. For $\vec{w}=(w_1,\ldots,w_m)$, set
$$
\nu_{\vec{w}}=\prod_{j=1}^m w_j^{p/p_j}.
$$
For $\theta\geq 0$, we say that $\vec{w}$ is in the class $A_{\vec{p}}^\theta$ if
$$
\Big(\f{1}{|B|}\int_B \nu_{\vec{w}}(x)dx\Big)^{1/p}\prod_{j=1}^m \Big(\f{1}{|B|}\int_B w_j^{1-p_j'}(x)dx\Big)^{1/p'_j}\leq C(1+r_B)^\theta
$$
for all balls $B$. When $p_j=1$, $\Big(\f{1}{|B|}\int_B w_j^{1-p_j'}(x)dx\Big)^{1/p'_j}$ is understood $(\inf_{x\in Q}w_j(x))^{-1}$.
\end{defn}

For $1\leq p_1, \ldots, p_m<\vc$, we set $A^\vc_{\vec{p}}=\cup_{\theta\geq 0}A^\theta_{\vec{p}}$.

\medskip

When $\theta=0$, the class $A^0_{\vec{p}}$ coincides with the class of multiple weights $A_{\vec{p}}$ introduced by \cite{LOPTG}. The following result gives a characterization of the class  $A^\vc_{\vec{p}}$ whose proof is similar to that of \cite[Theorem 3.6]{LOPTG}.
\begin{prop}\label{prop1 A_p}
Let $1\leq p_1, \ldots, p_m<\vc$ and $\vec{w}=(w_1, \ldots, w_m)$. Then the following statements are equivalent:

(i) $\vec{w}\in A^\vc_{\vec{p}}$;

(ii) $w_j^{1-p_j'}\in A^\vc_{mp_j'}, j=1,\ldots, m$ and $\nu_{\vec{w}}\in A^\vc_{mp}$.

\end{prop}

Note that the class $A_{\vec{p}}^\vc$ is not increasing. It means that for $\vec{p}=(p_1, \ldots, p_m)$ and $\vec{q}=(q_1, \ldots, q_m)$ with $p_j\leq q_j, j=1,\ldots, m$, the following may not be true $A^\vc_{\vec{p}}\subset A^\vc_{\vec{q}}$, see \cite[Remark 7.3]{LOPTG}. However, we have the following result.

\begin{prop}\label{prop2 A_p}
Let $1\leq p_1, \ldots, p_m<\vc$ and $\vec{w}=(w_1, \ldots, w_m) \in A^\vc_{\vec{p}}$. Then,

(i) For any $r\geq 1$, $\vec{w}\in A^\vc_{r\vec{p}}$;

(ii) If $1< p_1, \ldots, p_m<\vc$, then there exists $r>1$ so that $\vec{w}\in A^\vc_{\vec{p}/r}$.
\end{prop}

\emph{Proof:} (i) Assume that $\vec{w}\in A^\theta_{\vec{p}}$ for some $\theta\geq 0$. By definition, there exists $C>0$ so that for all balls $B$ there holds,
$$
\Big(\f{1}{|B|}\int_B \nu_{\vec{w}}(x)dx\Big)^{1/p}\prod_{j=1}^m \Big(\f{1}{|B|}\int_B w_j^{-\f{1}{p_j-1}}(x)dx\Big)^{\f{p_j-1}{p_j}}\leq C(1+r_B)^\theta.
$$
For $r>1$, by H\"older inequality we have
$$
\Big(\f{1}{|B|}\int_B w_j^{-\f{1}{rp_j-1}}(x)dx\Big)^{\f{rp_j-1}{rp_j}}\leq \Big(\f{1}{|B|}\int_B w_j^{-\f{1}{p_j-1}}(x)dx\Big)^{\f{p_j-1}{rp_j}}.
$$
Therefore, for all balls $B$, we have
$$
\Big(\f{1}{|B|}\int_B \nu_{\vec{w}}(x)dx\Big)^{1/rp}\prod_{j=1}^m \Big(\f{1}{|B|}\int_B w_j^{-\f{1}{rp_j-1}}(x)dx\Big)^{\f{rp_j-1}{rp_j}}\leq C(1+r_B)^{\theta/r}.
$$
This implies that $\vec{w}\in A^\vc_{r\vec{p}}$.

\vskip0.75cm

(ii) We exploit some ideas in \cite{LOPTG} to our situation. In the light of Propositions \ref{property-Avc} and \ref{prop1 A_p}, we can pick $r_0>1$ and $\eta\geq 0$ so that
\begin{equation}\label{eq1-A_pproperty}
\Big(\f{1}{|B|}\int_B w_j^{-\f{r_0}{p_j-1}}(x)\Big)^{1/r_0}\leq C\Big(\f{1}{|B|}\int_B w_j^{-\f{1}{p_j-1}}(x)\Big)(1+r_B)^\theta
\end{equation}
for all balls $B$ and $j=1,\ldots, m$.

Taking $r>1$ so that $r<\f{r_0p_j}{p_j+r_0-1}$ for all $j=1,\ldots,m$, then we have for all $j$,
$$
\f{r_0(p_j-r)}{r(p_j-1)}<1.
$$
This together with H\"older inequality gives
$$
\Big(\f{1}{|B|}\int_B w_j^{-\f{r}{p_j-r}}(x)dx\Big)^{\f{p_j-r}{p_j}}\leq \Big(\f{1}{|B|}\int_B w_j^{-\f{r_0}{p_j-1}}(x)dx\Big)^{\f{r(p_j-1)}{r_0p_j}}.
$$
Due to (\ref{eq1-A_pproperty}), we have
$$
\Big(\f{1}{|B|}\int_B w_j^{-\f{r}{p_j-r}}(x)dx\Big)^{\f{p_j-r}{p_j}}\leq \Big(\f{1}{|B|}\int_B w_j^{-\f{1}{p_j-1}}(x)dx\Big)^{\f{r(p_j-1)}{p_j}}(1+r_B)^{r\theta}.
$$
Therefore,
$$
\Big(\f{1}{|B|}\int_B \nu_{\vec{w}}(x)dx\Big)^{r/p}\prod_{j=1}^m \Big(\f{1}{|B|}\int_B w_j^{-\f{1}{p_j/r-1}}(x)dx\Big)^{\f{p_j/r-1}{p_j/r}}\leq C(1+r_B)^{r\theta}.
$$
It yields $\vec{w}\in A^\vc_{\vec{p}/r}$.
\begin{flushright}
    $\Box$
\end{flushright}

\subsubsection{New BMO function spaces BMO$_\vc$}

In this section, we will recall the definition and some basic properties of the new BMO function spaces. According to \cite{BHS1}, the new BMO space $BMO_\theta$ with $\theta\geq 0$ is defined as a set of all locally integrable functions $b$ satisfying
\begin{equation}\label{eq1-intro}
\f{1}{|B|}\int_B|b(y)-b_B|dy\leq C(1+r_B)^\theta
\end{equation}
where $B=B(x_B,r_B)$ and $b_B=\f{1}{|B|}\int_B b$. A norm for $b\in BMO_\theta$, denoted by $\|b\|_\theta$, is given by the infimum of the constants satisfying (\ref{eq1-intro}). Clearly $BMO_{\theta_1}\subset BMO_{\theta_2}$ for $\theta_1\leq \theta_2$ and $BMO_0=BMO$. We define $BMO_\vc=\cup_{\theta\geq 0}BMO_\theta$.

The following result can be considered to be a variant of John-Nirenberg inequality for the spaces $BMO_\vc$.
\begin{prop}\label{JNforBMOL}
Let $\theta>0, s\geq 1$. If $b\in BMO_\theta$ then for all $B=(x_0, r)$

i) $$
\Big(\f{1}{|B|}\int_B|b(y)-b_B|^sdx\Big)^{1/s}\lesssim \|b\|_{\theta}(1+r_B)^{\theta};
$$

ii) $$
\Big(\f{1}{|2^kB|}\int_{2^kB}|b(y)-b_B|dx\Big)^{1/s}\lesssim \|b\|_{\theta} k(1+2^kr_B)^{\theta}
$$
for all $k\in \mathbb{N}$.

\end{prop}

We refer to Proposition 3 in \cite{BHS2} for the proof.

\subsection{Weighted estimates for some maximal operators}

A ball of the form $B(x_B,r_B)$ is called {\it a critical ball} if $r_B =1$. We have the following result.

\begin{prop}\label{coveringlemma}
There exists a sequence of points $x_j, j\geq 1$ in $\mathbb{R}^n$ so that the family of critical balls $\{Q_j\}_j$ where $Q_j:=B(x_j, 1)$, $j\geq 1$ satisfies

(i) $\cup_j Q_j = \mathbb{R}^n$.

(ii) There exists a constant $C$ such that for any $\sigma>1$, $\sum_j \chi_{\sigma Q_j}\leq C\sigma^{n}$.
\end{prop}

For the proof, we refer the reader to  \cite{B} (see also \cite{DZ}).

\medskip

We consider the following maximal functions for $g\in L^1_{{\rm loc}}(\mathbb{R}^n)$ and $x\in \mathbb{R}^n$
$$
M_{{\rm loc}, \alpha}g(x)=\sup_{x\in B\in \mathcal{B}_{\alpha}}\f{1}{|B|}\int_B|g|,
$$
$$
M^\sharp_{{\rm loc}, \alpha}g(x)=\sup_{x\in B\in \mathcal{B}_{\alpha}}\f{1}{|B|}\int_B|g-g_B|\approx \inf_{c\in \mathbb{R}}\sup_{x\in B\in \mathcal{B}_{\alpha}}\f{1}{|B|}\int_B|g-c|,
$$
where $\mathcal{B}_{\alpha}=\{B(y,r): y\in \mathbb{R}^n \ \text{and} \ r\leq \alpha \}$.

Also, given a ball $Q$, we define the following maximal functions for $g\in L^1_{{\rm loc}}(\mathbb{R}^n)$ and $x\in Q$
$$
M_{Q}g(x)=\sup_{x\in B\in \mathcal{F}(Q)}\f{1}{|B\cap Q|}\int_{B\cap Q}|g|,
$$
$$
M^\sharp_{Q}g(x)=\sup_{x\in B\in \mathcal{F}(Q)}\f{1}{|B\cap Q|}\int_{B\cap Q}|g-g_{B\cap Q}|,
$$
where $\mathcal{F}(Q)=\{B(y,r): y\in Q, r>0 \}$.

We have the following lemma.
\begin{lem}\label{FSinequalityversion}
For $0< p<\vc$, then there exists $\beta$ such that if $\{Q_k\}_{k}$ is a sequence of balls as in Proposition \ref{coveringlemma} then for all $g\in L^1_{{\rm loc}}(\mathbb{R}^n)$ and $w\in A^\vc_{\vc}$, we have

(i)$$
\int_{\mathbb{R}^n}|M_{{\rm loc},\beta}g(x)|^pw(x)dx \lesssim \int_{\mathbb{R}^n}|M^\sharp_{{\rm loc},4}g(x)|^pw(x)dx+\sum_{k}w(2Q_k)\Big(\f{1}{|2Q_k|}\int_{2Q_k}|g|\Big)^p;
$$
and

(ii) $\|M_{{\rm loc},\beta}g\|^p_{L^{p,\vc}(w)}\leq \|M^\sharp_{{\rm loc},4}g\|^p_{L^{p,\vc}(w)}+\sum_{k}w(2Q_k)\Big(\f{1}{|2Q_k|}\int_{2Q_k}|g|\Big)^p$.
\end{lem}
\emph{Proof:} We refer to \cite[Lemma 2.4]{B} for the proof of (i). The proof of (ii) is similar to that of (i) and we omit details here.
\begin{flushright}
    $\Box$
\end{flushright}

\medskip

Throughout this paper, we always assume that $N$ is a sufficiently large number and different from line to line. For $\kappa\geq 1$, $p> 0$,  $\vec{f}=(f_1,\ldots, f_m), f_j\in L^1_{{\rm loc}}(\mathbb{R}^n)$ for all $j=1,\ldots, m$ and $x\in \mathbb{R}^n$, we define the maximal function $\mathcal{M}$ by
$$
\mathfrak{M}_{\kappa, p}(\vec{f})(x)=\sup_{Q\ni x; Q \ {\rm is \ critical}}\sum_{k=0}^\vc 2^{-Nk}\prod_{j=1}^m \Big(\f{1}{|2^k\widehat{Q}|}\int_{2^k\widehat{Q}}|f_j(z)|^pdz\Big)^{1/p}
$$
where $\widehat{Q}=\kappa Q$.

For simplicity, we shall write $\mathfrak{M}_p$ and $\mathfrak{M}$ instead of $\mathfrak{M}_{1,p}$ and $\mathfrak{M}_{1,1}$, respectively. The following result gives the  weighted estimates for $\mathfrak{M}_{\kappa,p}$.
\begin{prop}\label{weighted estiamtes for G}
Let $p_1, \ldots, p_m \geq s >0$ and $w\in A_{\vec{p}/s}^\vc$. Then we have
$$
\|\mathfrak{M}_{\kappa, s}(\vec{f})\|_{L^p(\nu_{\vec{w}})}\lesssim \prod_{j=1}^m \|f\|_{L^{p_j}(w_j)}.
$$
\end{prop}
\emph{Proof:} Without of the loss of generality, we assume that $\alpha =1$. Let $\{Q_\ell\}$ be the family of balls as in Proposition \ref{coveringlemma}. Then we have
$$
\begin{aligned}
\|\mathfrak{M}_{s}(\vec{f})\|_{L^p(\nu_{\vec{w}})}\leq C \Big(\sum_\ell \int_{Q_\ell}|\mathfrak{M}_{s}(\vec{f})(x)|^p \nu_{\vec{w}}dx\Big)^{1/p}.
\end{aligned}
$$
If $x\in Q_\ell\cap Q$, where $Q$ is a critical ball, then $2^kQ\subset 2^{k+1}Q_\ell$ and $|2^kQ|\approx |2^{k+1}Q_\ell|$. Therefore, for $x\in Q_\ell$ we have
$$
\mathfrak{M}_{s}(\vec{f})(x)\leq C\sum_{k=0}^\vc 2^{-Nk}\prod_{j=1}^m \Big(\f{1}{|2^kQ_\ell|}\int_{2^kQ_\ell}|f_j(z)|^sdz\Big)^{1/s}.
$$
So,
$$
\begin{aligned}
\|\mathfrak{M}_{s}(\vec{f})\|_{L^p(\nu_{\vec{w}})}&\leq C \Big(\sum_\ell \int_{Q_\ell}\Big| \sum_{k=0}^\vc 2^{-Nk}\prod_{j=1}^m \Big(\f{1}{|2^kQ_\ell|}\int_{2^kQ_\ell}|f_j(z)|^sdz\Big)^{1/s} \Big|^p \nu_{\vec{w}}dx\Big)^{1/p}\\
&\leq C \sum_{k=0}^\vc 2^{-Nk}\Big(\sum_\ell \nu_{\vec{w}}(Q_\ell)\prod_{j=1}^m \Big(\f{1}{|2^kQ_j|}\int_{2^kQ_j}|f_j(z)|^sdz\Big)^{p/s}\Big)^{1/p}\\
&\leq C \sum_{k=0}^\vc 2^{-Nk}\Big(\sum_\ell \nu_{\vec{w}}(2^kQ_\ell)\prod_{j=1}^m \Big(\f{1}{|2^kQ_j|}\int_{2^kQ_j}|f_j(z)|^sdz\Big)^{p/s}\Big)^{1/p}.
\end{aligned}
$$
Assume that $\vec{w}\in A^\theta_{\vec{p}/s}$ for some $\theta\geq 0$.  For each $j$, by H\"older inequality, we have
$$
\begin{aligned}
\int_{2^kQ_\ell}|f_j(z)|^sdz\leq \|f_j\|^s_{L^{p_j}(w_j, 2^kQ_\ell)}\Big(\int_{2^kQ_\ell}w^{-\f{1}{p_j/s-1}}\Big)^{\f{p_j/s}{p_j/s -1}}.
\end{aligned}
$$
This together with definition of $ A^\theta_{\vec{p}}$ class gives
$$
\begin{aligned}
\|\mathfrak{M}_{s}(\vec{f})\|_{L^p(\nu_{\vec{w}})}&\leq C \sum_{k=0}^\vc 2^{-k(N-\theta/s)}\Big(\sum_\ell \prod_{j=1}^m \|f_j\|^p_{L^{p_j}(w_j, 2^kQ_\ell)}\Big)^{1/p}.
\end{aligned}
$$
Using H\"older inequality and (ii) of Proposition \ref{coveringlemma}, we get that
$$
\begin{aligned}
\sum_\ell \prod_{j=1}^m \|f_j\|^p_{L^{p_j}(w_j, 2^kQ_\ell)}&\leq \prod_{j=1}^m\Big(\sum_{\ell}\|f_j\|^{p_j}_{L^{p_j}(w_j, 2^kQ_\ell)}\Big)^{p/p_j}\\
&\leq C2^{-kn}\prod_{j=1}^m\|f_j\|^{p}_{L^{p_j}(w_j)}.
\end{aligned}
$$
This completes our proof.
\begin{flushright}
    $\Box$
\end{flushright}

\medskip

For a family of balls $\{Q_k\}_k$ given by Proposition \ref{coveringlemma}, for $s>0$, we define the operator $\mathcal{M}_{{\rm loc},s}$ by setting
\begin{equation}\label{defnoftidleM}
\mathcal{M}_{{\rm loc},s}(\vec{f})=\sum_k\chi_{Q_k}\mathcal{M}_s(\vec{f}\chi_{\widetilde{Q}_k})
\end{equation}
where $\widetilde{Q}_j= 8 Q_j, \vec{f}\chi_{\widetilde{Q}_k}=(f_1\chi_{\widetilde{Q}_k}, \ldots, f_m\chi_{\widetilde{Q}_k})$, and the maximal function $\mathcal{M}_s$ is defined by
$$
\mathcal{M}_s(\vec{f})(x)=\sup_{Q\ni x}\prod_{j=1}^m\Big(\f{1}{|Q|}\int_Q|f(z)|^sdz\Big)^{1/s}.
$$
When $s=1$, we drop the subindex $s$ to write $\mathcal{M}_{{\rm loc}}$ instead of $\mathcal{M}_{{\rm loc},1}$.

\medskip

We have the following result.
\begin{prop}\label{rem1}
(i) If $p_1, \ldots, p_m>s> 0$ and $\vec{w}\in A^\vc_{\vec{p}/s}$, then we have
$$
\|\mathcal{M}_{{\rm loc},s}(\vec{f})\|_{L^p(\nu_{\vec{w}})}\leq C\prod_{j=1}^m\|f_j\|_{L^{p_j}(w_j)}.
$$

(ii) If $p_1, \ldots, p_m\geq s>0$ and $\vec{w}\in A^\vc_{\vec{p}/s}$, then we have
$$
\|\mathcal{M}_{{\rm loc},s}(\vec{f})\|_{L^{p,\vc}(\nu_{\vec{w}})}\leq C\prod_{j=1}^m\|f_j\|_{L^{p_j}(w_j)}.
$$
\end{prop}
\emph{Proof:} The give the proof of (i) only. The proof of (ii) can be dealt by the analogous argument.

\medskip

We first have
$$
\begin{aligned}
\int_{\RR^n}|\mathcal{M}_{{\rm loc},s}(\vec{f})(x)|^p\nu_{\vec{w}}(x)dx&=\sum_{\ell}\int_{Q_\ell}|\mathcal{M}_s(\vec{f}\chi_{\widetilde{Q}_\ell})|^p\nu_{\vec{w}}(x)dx\\
&=\sum_{\ell}\int_{\RR^n}|\mathcal{M}_s(\vec{f}\chi_{\widetilde{Q}_\ell})|^p\nu_{\vec{w}_\ell}(x)dx
\end{aligned}
$$
where $\vec{w}_j = (w_1\chi_{Q_j}, \ldots, w_m\chi_{Q_j})$.

For each $j$, it can be verified that $\vec{w}_j \in A^0_{\vec{p}/s}$. Therefore, by Theorem 3.7 in \cite{LOPTG}, (ii) of Proposition \ref{coveringlemma} and H\"older inequality, we have
$$
\begin{aligned}
\int_{\RR^n}|\mathcal{M}_{{\rm loc},s}(\vec{f})(x)|^p\nu_{\vec{w}}(x)dx&\leq C\sum_{\ell}\prod_{j=1}^m\|f_j\|^p_{L^{p_j}(w_j, \widetilde{Q}_\ell)}\\
&\leq C\prod_{j=1}^m\Big(\sum_{\ell}\|f_j\|^{p_j}_{L^{p_j}(w_j, \widetilde{Q}_\ell)}\Big)^{p/p_j}\\
&\leq C\prod_{j=1}^m\|f_j\|^{p}_{L^{p_j}(w_j)}.
\end{aligned}
$$
This completes our proof.
\begin{flushright}
    $\Box$
\end{flushright}

When $m=1$, we write $\widetilde{M}_{{\rm loc}, s}$ instead of $\mathcal{M}_{{\rm loc}, s}$. As a direct consequence of Proposition \ref{rem1}, we have the following result.
\begin{prop}\label{rem2}
(i) If $p>s> 0$ and $w\in A^\vc_{p/s}$, then we have
$$
\|\widetilde{M}_{{\rm loc},s}(f)\|_{L^p(w)}\leq C\|f\|_{L^{p}(w)}.
$$

(ii) If $p=s$ and $w\in A^\vc_{1}$, then we have
$$
\|\widetilde{M}_{{\rm loc},s}(\vec{f})\|_{L^{s,\vc}(w)}\leq C\|f\|_{L^{s}(w)}.
$$
\end{prop}
\section{Main results}
\subsection{Weighted norm inequalities for multilinear operators}
In this section, we establish the weighted norm inequalities for multilinear operators in the product of weighted $L^p$ spaces by using the new class of weights $A^\vc_{\vec{p}}$.
\begin{thm}\label{thm1}
Let $T$ satisfy (H1), (H2) and (H3). Then the following statements hold:

(i) If $1<p_1,\ldots, p_m<\vc$ and $\vec{w}\in A_{\vec{p}}^\vc$, then
$$
\|T(\vec{f})\|_{L^p(\nu_{\vec{w}})}\leq C\prod_{j=1}^m \|f_j\|_{L^{p_j}(w_j)}.
$$

(ii) If $1\leq p_1,\ldots, p_m<\vc$ and at least one of the $p_j=1$, then
$$
\|T(\vec{f})\|_{L^{p,\vc}(\nu_{\vec{w}})}\leq C\prod_{j=1}^m \|f_j\|_{L^{p_j}(w_j)}.
$$
\end{thm}

To prove Theorem \ref{thm1}, we need the following auxiliary propositions.

\begin{prop}\label{pro1}
Let $T$ satisfy (H1), (H2) and (H3). Then for any critical ball $Q$ and $0<\delta<1/m$, we have
$$
\Big(\f{1}{|2Q|}\int_{2Q}|T(\vec{f})(z)|^\delta dz\Big)^{1/\delta}\leq C\inf_{y\in 2Q}\mathfrak{M}(\vec{f})(y).
$$
\end{prop}
\emph{Proof:} We split $f_j=f^0_j + f^\vc_j$ where $f_j^0=f_j\chi_{4Q}$ for all $j$. Then,
$$
T(\vec{f})(x)=T(\vec{f}_0)(x)+\sum_{\alpha\in \mathcal{I}}T(f_1^{\alpha_1},\ldots,f_m^{\alpha_m})(x)
$$
where $\vec{f}_0=(f^0_1,\ldots, f_m^0)$, $\alpha_1,\ldots, \alpha_m\in \{0,\vc\}$ and $\mathcal{I}:=\{\alpha=(\alpha_1,\ldots, \alpha_m): \alpha\neq (0,\ldots,0)\}$.

Therefore,
$$
\begin{aligned}
\Big(\f{1}{|2Q|}\int_{2Q}|T(\vec{f})(z)|^\delta dz\Big)^{1/\delta}&\leq C\Big(\f{1}{|2Q|}\int_{2Q}|T(\vec{f}_0)(z)|^\delta dz\Big)^{1/\delta}\\
&~~~~~+C\Big(\f{1}{|2Q|}\int_{2Q}\Big|\sum_{\alpha\in \mathcal{I}}T(f_1^{\alpha_1},\ldots,f_m^{\alpha_m})(z)\Big|^\delta dz\Big)^{1/\delta}\\
&:= I_1+I_2.
\end{aligned}
$$

Note that
$$
\|T(\vec{f}_0)\|_{L^{1/m,\vc}}\leq C\prod_{j=1}^m \|f_j^0\|_{L^{1}}.
$$
This estimate together with Kolmogorov inequality implies
$$
\begin{aligned}
I_1&\leq \|T(\vec{f}_0)\|_{L^{1/m,\vc}(2Q, \f{dz}{|2Q|})} \leq C\inf_{y\in 2Q}\mathfrak{M}(\vec{f})(y).
\end{aligned}
$$

\medskip

Let us take care $I_2$. For each $\alpha=(\alpha_1, \ldots, \alpha_m)\in \mathcal{I}$ and $z\in 2Q$, we have, by (H1),
$$
\begin{aligned}
|T(f_1^{\alpha_1},\ldots,f_m^{\alpha_m})(z)|&\leq C\int_{\RR^{mn}}\f{|f_1^{\alpha_1}(y_1)\ldots f_m^{\alpha_m}(y_m)|}{(|z-y_1|+\ldots+|z-y_m|)^{mn+N}}dy\\
&\leq C\sum_{k\geq 2} \int_{(2^{k+1}Q)^m\backslash (2^{k}Q)^m}\f{|f_1(y_1)\ldots f_m(y_m)|}{(|z-y_1|+\ldots+|z-y_m|)^{mn+N}}dy\\
&\leq C\sum_{k\geq 2} 2^{-kN} \f{1}{|2^{k+1}Q|^m}\int_{(2^{k+1}Q)^m}|f_1(y_1)|\ldots |f_m(y_m)|dy\\
&\leq C\inf_{y\in 2Q} \mathfrak{M}(\vec{f})(y).
\end{aligned}
$$
This completes our proof.
\begin{flushright}
    $\Box$
\end{flushright}

\medskip

\begin{prop}\label{pro2}
Let $T$ satisfy (H1), (H2) and (H3), let $0<\delta<1/m$. Then we have
$$
(M^\sharp_{{\rm loc}, 4}|T(\vec{f}|^{\delta})(x))^{1/\delta}\leq C\mathfrak{M}(\vec{f})(x)+C\mathcal{M}_{{\rm loc}}(\vec{f})(x).
$$
\end{prop}
\emph{Proof:} It suffices to show that for each ball $B\ni x$ with $r_B\leq 4$, we can pick a number $c_B$ so that
$$
\Big(\f{1}{|B|}\int_B\Big||T(\vec{f})(z)|^{\delta} -|c_B|^\delta\Big| dz\Big)^{1/\delta}\leq C\mathfrak{M}(\vec{f})(x)+C\mathcal{M}_{{\rm loc}}(\vec{f})(x).
$$
By the same decomposition as in Proposition \ref{pro1}, we can write
$$
\begin{aligned}
\Big(\f{1}{|B|}\int_B\Big||T(\vec{f})(z)|^{\delta} -|c_B|^\delta\Big| dz\Big)^{1/\delta}&\leq C \Big(\f{1}{|B|}\int_B|T(\vec{f})(z) -c_B|^\delta dz\Big)^{1/\delta}\\
&\leq C\Big(\f{1}{|B|}\int_B|T(\vec{f}_0)(z)|^{\delta}dz\Big)^{1/\delta}\\
&~~+C\Big(\f{1}{|B|}\int_B\Big|\sum_{\alpha\in \mathcal{I}}T(f_1^{\alpha_1}, \ldots, f_m^{\alpha_m})(z)-c_B\Big|^{\delta}dz\Big)^{1/\delta}\\
&=E_1+E_2.
\end{aligned}
$$
By Kolmogorov inequality and the fact that $T$ is bounded form $L^1\times\ldots\times L^1$ into $L^{1/m,\vc}$, we have
$$
E_1\leq C\prod_{j=1}^m \f{1}{|4B|}\int_{4B}|f_j(z)|dz.
$$
Let $\{Q_\ell\}_{\ell}$ be the family in Proposition \ref{coveringlemma}. Then if $4B\cap Q_{\ell}\neq \emptyset$, we have $4B\subset \widetilde{Q}_\ell: = 8Q_\ell$. Therefore,
$$
E_1\leq C \mathcal{M}_{{\rm loc}}(\vec{f})(x).
$$

For the second term $E_2$. Taking $c_B=\sum_{\alpha\in \mathcal{I}}T(f_1^{\alpha_1}, \ldots, f_m^{\alpha_m})(x_B)$, we have
$$
\sum_{\alpha\in \mathcal{I}}T(f_1^{\alpha_1}, \ldots, f_m^{\alpha_m})(z)-c_B=\sum_{\alpha\in \mathcal{I}}\Big[T(f_1^{\alpha_1}, \ldots, f_m^{\alpha_m})(z)-T(f_1^{\alpha_1}, \ldots, f_m^{\alpha_m})(x_B)\Big].
$$
For each $\alpha\in \mathcal{I}$ and $N>0$, by (H2), we write
$$
\begin{aligned}
\Big|T(f_1^{\alpha_1}, &\ldots, f_m^{\alpha_m})(z)-T(f_1^{\alpha_1}, \ldots, f_m^{\alpha_m})(x_B)\Big|\\
&\leq C\int_{\RR^{mn}}|K(z,y_1,\ldots,y_m)-K(x_B,y_1,\ldots,y_m)|\prod_{j=1}^m |f^{\alpha_j}_j(y_j)|dy\\
&\leq C\sum_{k\geq 2}\int_{(2^{k+1}B)^m\backslash (2^{k}B)^m}|K(z,y_1,\ldots,y_m)-K(x_B,y_1,\ldots,y_m)|\prod_{j=1}^m |f_j(y_j)|dy\\
&\leq C\sum_{k\geq 2}\min\{1, (2^kr_B)^{-N}\} \f{2^{-k\epsilon}}{|2^{k+1}B|^m}\int_{(2^{k+1}B)^m}\prod_{j=1}^m |f_j(y_j)|dy\\
&\leq C\sum_{k\geq k_0}\ldots + C\sum_{k< k_0}\ldots:= E_{21}+E_{22}
\end{aligned}
$$
where $k_0$ is the smallest integer so that $2^{k_0}r_B\geq 1$.

Let us estimate $E_{22}$ first. Since $r_{2^{k+1}B}\leq 4$ for all $k<k_0$, the similar argument used in the estimate $E_1$ gives
$$
\f{1}{|2^{k+1}B|^m}\int_{(2^{k+1}B)^m}\prod_{j=1}^m |f_j(y_j)|dy\leq C\mathcal{M}_{{\rm loc}}(\vec{f})(x)
$$
for all $k<k_0$.

Therefore,
$$
E_{22}\leq C\sum_{k< k_0}2^{-k\epsilon}\mathcal{M}_{{\rm loc}}(\vec{f})(x)\leq C\mathcal{M}_{{\rm loc}}(\vec{f})(x).
$$

For the term $E_{21}$, setting $\widehat{B}=2^{k_0}B$, then $1\leq r_{\widehat{B}}<2$. We have
$$
\begin{aligned}
E_{21}&= \sum_{k\geq k_0}(2^{k-k_0}2^{k_0}r_B)^{-N}\f{2^{-k\epsilon}}{|2^{k-k_0+1}\widehat{B}|^m}\int_{(2^{k-k_0+1}\widehat{B})^m}\prod_{j=1}^m |f_j(y_j)|dy\\
&= \sum_{k\geq k_0}2^{-kN}\f{1}{|2^{k}\widehat{B}|^m}\int_{(2^{k}\widehat{B})^m}\prod_{j=1}^m |f_j(y_j)|dy\\
&\leq C\mathfrak{M}(\vec{f})(x).
\end{aligned}
$$
This completes our proof.
\begin{flushright}
    $\Box$
\end{flushright}

\vskip1cm

We are now in position to prove Theorem \ref{thm1}.

\vskip1cm

\emph{Proof of Theorem \ref{thm1}:} Since the proofs of (i) and (ii), we give only the proof of (i).

(i) Since $\nu_{\vec{w}}\in A^\vc_{mp}$, using Proposition \ref{FSinequalityversion} for $\delta<p<1/m$, we have
$$
\begin{aligned}
\|T(\vec{f})\|_{L^p(\nu_{\vec{w}})}^p&:=\|(T(\vec{f}))^\delta\|_{L^{p/\delta}(\nu_{\vec{w}})}^{p/\delta}\\
&\leq  \int_{\RR^n}|M_{{\rm loc},\beta}(|T(\vec{f})|^\delta)(x)|^{p\delta}\nu_{\vec{w}}(x)dx\\
&\leq C\int_{\RR^n}|M^\sharp_{{\rm loc},4}(|T(\vec{f})|^\delta)(x)|^{p\delta}\nu_{\vec{w}}(x)dx+C\sum_{\ell}\nu_{\vec{w}}(2Q_\ell)\Big(\f{1}{2Q_k}\int_{2Q_k}|T(\vec{f})(z)|^{\delta}dz\Big)^{p/\delta}\\
&= E_1+E_2.
\end{aligned}
$$

By Proposition \ref{pro1}, Proposition \ref{coveringlemma} and Proposition \ref{weighted estiamtes for G}, we have
$$
\begin{aligned}
E_2&\leq C \sum_{\ell}\nu_{\vec{w}}(2Q_\ell)\Big[\inf_{y\in 2Q_\ell}\mathfrak{M}(\vec{f})(y)\Big]^p\\
&\leq C \sum_{\ell}\int_{Q_\ell}|\mathfrak{M}(\vec{f})(z)|^p\nu_{\vec{w}}(z)dz\\
&\leq C\|\mathfrak{M}(\vec{f})\|^p_{L^p(\nu_{\vec{w}})}\\
&\leq C\prod_{j=1}^m \|f_j\|^p_{L^{p_j}(w_j)}.
\end{aligned}
$$

It remains to estimate $E_1$. In the light of Proposition \ref{pro2}, Proposition \ref{weighted estiamtes for G} and Proposition \ref{rem1}, we have
$$
\begin{aligned}
E_1&\leq C \int_{\RR^n}|\mathfrak{M}(\vec{f})(z)|^p\nu_{\vec{w}}(z)dz + C\int_{\RR^n}|\mathcal{M}_{{\rm loc}}(\vec{f})(z)|^p\nu_{\vec{w}}(z)dz\\
&\leq C\prod_{j=1}^m \|f_j\|^p_{L^{p_j}(w_j)}.
\end{aligned}
$$
This completes our proof.
\begin{flushright}
    $\Box$
\end{flushright}

\subsection{Weighted norm inequalities for the commutators of multilinear operators with BMO$_\vc$ functions}

Let $\vec{b}=(b_1, \ldots, b_m)$, where $b_j$ is a locally integrable function for  $j=1,\ldots,m$. We consider the $m$-linear commutator of $T$ and $\vec{b}$
$$
T_{\vec{b}}(\vec{f})=\sum_{j=1}^m T^j_{\vec{b}}(\vec{f})
$$
where
$$
T^j_{\vec{b}}(\vec{f})=b_jT(\vec{f})-T(f_1,\ldots, b_jf_j,\ldots, f_m).
$$

Let $\vec{\theta}=(\theta_1,\ldots, \theta_m)$, $\theta_j\geq 0$ for all $j=1,\ldots,m$. For $\vec{b}\in BMO_{\vec{\theta}}$, we shall mean $b_j\in BMO_{\theta_j}$ for all $j=1,\ldots,m$ and $\|\vec{b}\|_{\vec{\theta}}=\sum_{j=1}^m \|b_j\|_{\theta_j}$.

Before coming to detail information, we need the following auxiliary results.

\begin{prop}\label{pro1-com}
Let $T$ satisfy (H1), (H2) and (H3) and $\vec{b}\in BMO_{\vec{\theta}}$. For any $p>1$, the following holds for any critical ball $Q$ and $0<\delta<1/m$,
$$
\Big(\f{1}{|2Q|}\int_{2Q}|T_{\vec{b}}(\vec{f})(z)|^\delta dz\Big)^{1/\delta}\leq C\|\vec{b}\|_{\vec{\theta}}\inf_{y\in 2Q}\mathfrak{M}_p(\vec{f})(y).
$$
\end{prop}
\emph{Proof:} By linearity it is suffices to consider the commutator with one symbol as follows:
$$
T_{b}(\vec{f})=bT(\vec{f})-T(bf_1,\ldots, f_m), b\in BMO_\theta.
$$
For any balls $Q$, we can write
$$
T_{b}(\vec{f})=(b-b_Q)T(\vec{f})-T((b-b_Q)f_1,\ldots, f_m):=I_1+I_2.
$$
Let $\delta<\delta'<1/m$ and $1/s+1/\delta'=1/\delta$ so that $s>1$. By H\"older inequality and Proposition \ref{pro1}, we have
$$
\begin{aligned}
I_1&\leq \Big(\f{1}{|2Q|}\int_{2Q}|T(\vec{f})(z)|^{\delta'} dz\Big)^{1/\delta'}\Big(\f{1}{|2Q|}\int_{2Q}|b(z)-b_Q|^s dz\Big)^{1/s}\\
&\leq C\|b\|_{\theta}\inf_{y\in 2Q}\mathfrak{M}(\vec{f})(y)\\
&\leq C\|b\|_{\theta}\inf_{y\in 2Q}\mathfrak{M}_p(\vec{f})(y).
\end{aligned}
$$

To estimate $I_2$, using the decomposition as in the proof of Proposition \ref{pro1}, we have
$$
T((b-b_Q)f_1, \ldots, f_m)(z)=T((b-b_Q)f^0_1, \ldots, f^0_m)(z)+\sum_{\alpha\in \mathcal{I}}T((b-b_Q)f_1^{\alpha_1},\ldots,f_m^{\alpha_m})(z)
$$
Using Kolmogorov inequality and the fact that $T$ maps continuously from $L^1\times \ldots L^1$ into $L^{1/m,\vc}$, we have
$$
\begin{aligned}
\Big(\f{1}{|2Q|}\int_{2Q}&|T((b-b_Q)f^0_1, \ldots, f^0_m)(z)|^\delta dz\Big)^{1/\delta}\\
&\leq \|T((b-b_Q)f^0_1, \ldots, f^0_m)\|_{L^{1/m,\vc}(2Q, \f{dz}{|2Q|})}\\
&\leq C\Big(\f{1}{|2Q|}\int_{4Q}|(b(z)-b_Q)f^0_1(z)|dz\Big)\prod_{j\neq 1}\Big(\f{1}{|2Q|}\int_{4Q}|f^0_j(z)|dz\Big)\\
&\leq C\Big(\f{1}{|2Q|}\int_{4Q}|(b(z)-b_Q)|^{p'}dz\Big)^{1/p'}\prod_{j}\Big(\f{1}{|2Q|}\int_{4Q}|f^0_j(z)|^pdz\Big)^{1/p}\\
&\leq C\|b\|_{\theta}\inf_{y\in 2Q}\mathfrak{M}_p(\vec{f})(y).
\end{aligned}
$$

For each $\alpha=(\alpha_1, \ldots, \alpha_m)\in \mathcal{I}$ and $z\in 2Q$, we have, by (H1), H\"older inequality and Proposition \ref{JNforBMOL},
$$
\begin{aligned}
|T((b-b_Q)&f_1^{\alpha_1},\ldots,f_m^{\alpha_m})(z)|\\
&\leq C\int_{\RR^{mn}}\f{|(b-b_Q)f_1^{\alpha_1}(y_1)\ldots f_m^{\alpha_m}(y_m)|}{(|z-y_1|+\ldots+|z-y_m|)^{mn+N}}dy\\
&\leq C\sum_{k\geq 2} \int_{(2^{k+1}Q)^m\backslash (2^{k}Q)^m}\f{|(b-b_Q)f_1(y_1)\ldots f_m(y_m)|}{(|z-y_1|+\ldots+|z-y_m|)^{mn+N}}dy\\
&\leq C\sum_{k\geq 2} 2^{-kN} \Big(\f{1}{|2^{k+1}Q|}\int_{2^{k+1}Q}|(b-b_Q)f_1(y_1)|dy_1\Big) \prod_{j\neq 1}\Big(\f{1}{|2^{k+1}Q|}\int_{2^{k+1}Q}|f_j(y_j)|dy_j\Big)\\
&\leq C\sum_{k\geq 2} 2^{-k(N-\theta)}\|b\|_\theta \Big(\f{1}{|2^{k+1}Q|}\int_{2^{k+1}Q}|f_1(y_1)|^pdy_1\Big)^{1/p} \prod_{j\neq 1}\Big(\f{1}{|2^{k+1}Q|}\int_{2^{k+1}Q}|f_j(y_j)|dy_j\Big)\\
&\leq C\|b\|_{\theta}\inf_{y\in 2Q} \mathfrak{M}_p(\vec{f})(y).
\end{aligned}
$$
This completes our proof.
\begin{flushright}
    $\Box$
\end{flushright}

\medskip

\begin{prop}\label{pro2-com}
Let $T$ satisfy (H1), (H2) and (H3), let $0<\delta<s<1/m$ and $\vec{b}\in BMO_{\vec{\theta}}$. Then we have, for any $p>1$,
$$
(M^\sharp_{{\rm loc}, 4}(|T_{\vec{b}}(\vec{f})|^{\delta})(x))^{1/\delta}\leq C\|\vec{b}\|_{\vec{\theta}}\Big[\widetilde{M}_{{\rm loc},s}(T(\vec{f}))(x)+\mathfrak{M}_p(\vec{f})(x)+\mathcal{M}_{{\rm loc},p}(\vec{f})(x)\Big]
$$
for all $x\in \RR^n$.
\end{prop}
\emph{Proof:} We need only to consider the commutator with one symbol as follows:
$$
T_{b}(\vec{f})=bT(\vec{f})-T(bf_1,\ldots, f_m)
$$
where $b\in BMO_\theta$ and $\lambda\in \RR$.

It suffices to show that for each ball $B\ni x$ with $r_B\leq 4$,
$$
\Big(\f{1}{|B|}\int_B\Big||T_b(\vec{f})(z)|^{\delta} -|c_B|^\delta\Big| dz\Big)^{1/\delta}\leq C\mathfrak{M}(\vec{f})(x)+C\mathcal{M}_{{\rm loc}}(\vec{f})(x)
$$
where $c_B$ is a constant which will be fixed later.

We write
$$
T_{b}(\vec{f})=(b-b_B)T(\vec{f})-T((b-b_B)f_1,\ldots, f_m)
$$
By the same decomposition as in Proposition \ref{pro1}, we can write
$$
\begin{aligned}
\Big(\f{1}{|B|}\int_B\Big||T_b(\vec{f})(z)|^{\delta} -|c_B|^\delta\Big| dz\Big)^{1/\delta}&\leq \Big(\f{1}{|B|}\int_B |T_b(\vec{f})(z)-c_B|^\delta dz\Big)^{1/\delta}\\
&\leq C\Big(\f{1}{|B|}\int_B|(b-b_B)T(\vec{f})(z)|^{\delta}dz\Big)^{1/\delta}\\
&~~+C\Big(\f{1}{|B|}\int_B|T((b-b_B)f_1^{0}, \ldots, f_m^{0})(z)|^{\delta}dz\Big)^{1/\delta}\\
&~~+C\Big(\f{1}{|B|}\int_B\Big|\sum_{\alpha\in \mathcal{I}}T((b-b_B)f_1^{\alpha_1}, \ldots, f_m^{\alpha_m})(z)-c_B\Big|^{\delta}dz\Big)^{1/\delta}\\
&=E_1+E_2+E_3.
\end{aligned}
$$
By H\"older inequality and Proposition \ref{JNforBMOL}, we have
$$
E_1\leq C\|b\|_{\theta}\Big(\f{1}{|4B|}\int_{B}|T(\vec{f})(z)|^sdz\Big)^{1/s}.
$$
Repeating the argument as in Proposition \ref{pro2}, we conclude that
$$
E_1\leq C\|b\|_{\theta}\widetilde{M}_{{\rm loc}, s}(T(\vec{f}))(x).
$$

The arguments in Propositions \ref{pro2} and \ref{pro1-com} yields that
$$
E_2\leq C\|b\|_{\theta}\prod_{j=1}^m \Big(\f{1}{|4B|}\int_{4B}|f_j(z)|^pdz\Big)^{1/p}\leq C\|b\|_{\theta}\mathcal{M}_{{\rm loc}, p}(\vec{f})(x).
$$

For the second term $E_3$. Taking $c_B=\sum_{\alpha\in \mathcal{I}}T((b-b_B)f_1^{\alpha_1}, \ldots, f_m^{\alpha_m})(x_B)$, we have
$$
\begin{aligned}
\sum_{\alpha\in \mathcal{I}}&T((b-b_B)f_1^{\alpha_1}, \ldots, f_m^{\alpha_m})(z)-c_B\\
&=\sum_{\alpha\in \mathcal{I}}\Big[T((b-b_B)f_1^{\alpha_1}, \ldots, f_m^{\alpha_m})(z)-T(f_1^{\alpha_1}, \ldots, f_m^{\alpha_m})(x_B)\Big].
\end{aligned}
$$
For each $\alpha\in \mathcal{I}$ and $N>0$, by (H2), we write
$$
\begin{aligned}
\Big|T((b-b_B)f_1^{\alpha_1}, &\ldots, f_m^{\alpha_m})(z)-T((b-b_B)f_1^{\alpha_1}, \ldots, f_m^{\alpha_m})(x_B)\Big|\\
&\leq C\int_{\RR^{mn}}|K(z,y_1,\ldots,y_m)-K(x_B,y_1,\ldots,y_m)||(b-b_B)f_1^{\alpha_1}(y_1)|\prod_{j\neq1}^m |f^{\alpha_j}_j(y_j)|dy\\
&\leq C\sum_{k\geq 2}\int_{(2^{k+1}B)^m\backslash (2^{k}B)^m}|K(z,y_1,\ldots,y_m)-K(x_B,y_1,\ldots,y_m)|\\
&~~~~~~~~~~~~~~\times|(b-b_B)f_1(y_1)|\prod_{j\neq1}^m |f_j(y_j)|dy\\
&\leq C\sum_{k\geq 2}\min\{1, (2^kr_B)^{-N}\} \f{2^{-k\epsilon}}{|2^{k+1}B|^m}\int_{(2^{k+1}B)^m}|(b-b_B)f_1(y_1)|\prod_{j\neq1}^m |f_j(y_j)|dy\\
&\leq C\sum_{k\geq k_0}\ldots + C\sum_{k< k_0}\ldots:= E_{31}+E_{32}
\end{aligned}
$$
where $k_0$ is the smallest integer so that $2^{k_0}r_B\geq 1$.

Let us estimate $E_{32}$ first. Since $r_{2^{k+1}B}\leq 4$ for all $k<k_0$, the similar argument used in the estimate $E_2$ gives
$$
\f{1}{|2^{k+1}B|^m}\int_{(2^{k+1}B)^m}|(b-b_B)f_1(y_1)|\prod_{j\neq1}^m |f_j(y_j)|dy\leq C\|b\|_{\theta}\mathcal{M}_{{\rm loc},p}(\vec{f})(x).
$$
Therefore,
$$
E_{32}\leq C\|b\|_{\theta}\sum_{k< k_0}2^{-k\epsilon}\mathcal{M}_{{\rm loc},p}(\vec{f})(x)\leq C\|b\|_{\theta}\mathcal{M}_{{\rm loc},p}(\vec{f})(x).
$$

For the term $E_{31}$, setting $\widehat{B}=2^{k_0}B$, then $1\leq r_{\widehat{B}}<2$. We have
$$
\begin{aligned}
E_{31}&= C\sum_{k\geq k_0}(2^{k-k_0}2^{k_0}r_B)^{-N}\f{2^{-k\epsilon}}{|2^{k-k_0+1}\widehat{B}|^m}\int_{(2^{k-k_0+1}\widehat{B})^m}|(b-b_B)f_1(y_1)|\prod_{j\neq1}^m |f_j(y_j)|dy\\
&= C\sum_{k\geq k_0}2^{-kN}\f{1}{|2^{k}\widehat{B}|^m}\int_{(2^{k}\widehat{B})^m}|(b-b_B)f_1(y_1)|\prod_{j\neq1}^m |f_j(y_j)|dy\\
&\leq C\sum_{k\geq k_0}2^{-k(N-\theta)}\Big(\f{1}{|2^{k}\widehat{B}|^m}\int_{(2^{k}\widehat{B})^m}\prod_{j=1}^m |f_j(y_j)|^pdy\Big)^{1/p}\\
&\leq C\|b\|_{\theta}\mathfrak{M}_p(\vec{f})(x).
\end{aligned}
$$
This completes our proof.
\begin{flushright}
    $\Box$
\end{flushright}

By similar arguments used in the proof of Theorem \ref{thm1}, we have the following result concerning the weighted norm inequality of the commutators $T_{\vec{b}}$.

\begin{thm}\label{thm2}
Let $T$ satisfy (H1), (H2) and (H3) and $\vec{b}\in BMO_{\vec{\theta}}$. Then we have, for  $1<p_1,\ldots, p_m<\vc$ and $\vec{w}\in A^\vc_{\vec{p}}$,
$$
\|T_{\vec{b}}(\vec{f})\|_{L^p(\nu_{\vec{w}})}\leq C\prod_{j=1}^m \|f_j\|_{L^{p_j}(w_j)}.
$$
\end{thm}

\section{Application to multilinear pseudodifferential operators}
In this section, we will apply the obtained results to study the weighted norm inequalities for multilinear pseudodifferential operators.

\medskip

Given a function $a: \RR^{n(m+1)}\rightarrow \mathbb{C}$ satisfying certain growth conditions, we define the multilinear operator $T_a$ to act on $m$ functions $f_1, \ldots, f_m\in \mathcal{S}(\RR^n)$ by setting
$$
T_a(\vec{f})(x)=\int_{\RR^{mn}}a(x,\xi)\prod_{j=1}^m \hat{f}_j(\xi_j)e^{i\langle x, \xi_1+\ldots+\xi_m \rangle}d\xi
$$
where $\xi=(\xi_1,\ldots, \xi_m)\in \RR^{mn}$.

We say that the symbol $a$ belongs to the H\"ormander class $m$-$S^l_{\rho,\delta}$, $l\in \RR, \rho, \delta\in [0,1]$ , i.e., for all multi-indices $\alpha, \beta_1, \ldots, \beta_m$ there holds
$$
|\partial_x^\alpha\partial_{\xi_1}^{\beta_1}\ldots \partial_{\xi_m}^{\alpha_m}a(x,\xi)|\leq C(1+|\xi_1|+\ldots+|\xi_1|)^{m+\delta |\alpha| -\rho(|\beta_1|+\ldots+|\beta_m|)}.
$$

We have the following result.
\begin{prop}\label{pro1-pseudo}
Let $a\in m$-$S^0_{1,\delta}$, $0\leq \delta<1$. Then $T_a$ satisfies (H1) and (H2).
\end{prop}
\emph{Proof:}
Let $\varphi_0:\RR^{mn} \rightarrow \RR$ be a smooth radial function which is equal to $1$ on the unit ball centered at origin and supported on its concentric double. Set $\varphi(\xi)=\varphi_0(\xi)-\varphi_0(2\xi)$ and $\varphi_k(\xi)=\varphi(2^{-k}\xi)$. Then, we have
$$
\sum_{k=0}^\vc\varphi_k(\xi)=1 \ \text{for all $\xi\in \RR^{mn}$}
$$
and supp $\varphi_k\subset \{\xi: 2^{k-1}\leq |\xi|\leq 2^{k+1}\}$ for all $k\geq 1$. Moreover, for  any multi-index $\alpha$ and $N\geq 0$, we have
\begin{equation}\label{eq1-pseudodiff}
|\pa^\alpha_\xi \varphi_k(\xi)|\leq c_\alpha 2^{-k|\alpha|}.
\end{equation}

Then we can write
$$
a(x,\xi)=\sum_{k=0}^\vc\varphi_k(\xi)a(x,\xi):=\sum_{k=0}^\vc a_k(x,\xi)
$$
It is not difficult to show that $T_a$ satisfies (H1). The proof of this part is standard and hence we omit details here.

It remains to check that $T_a$ satisfies (H1). To do this, we will work with each component $T_{a_k}$, for $k\geq 0$. With the same notations as in condition (H2), we consider two cases:

{\bf Case 1: $j\neq 0$}

In this situation, we can assume that $j=1$ and $\max_{k}|y_k-y_1|:=|y_0-y_1|$.

Let $K_k(y_0,y_1,\ldots, y_m)$ be the associated kernel of $T_{a_k}$. Then we have
$$
K_k(y_0,y_1,\ldots, y_m)=\int_{\RR^{mn}}a_k(y_0,\xi)\prod_{j=1}^me^{i\langle \xi_j, y_0-y_j \rangle}d\xi
$$
Therefore,
$$
\begin{aligned}
K_k(y_0,y_1,\ldots, y_m)&-K_k(y_0,y'_1,\ldots, y_m)\\
&=\int_{\RR^{mn}}a_k(y_0,\xi)\prod_{j\neq 1} e^{i\langle \xi_j, y_0-y_j \rangle}(e^{i\langle \xi_1, y_0-y_1 \rangle}-e^{i\langle \xi_1, y_0-y_1' \rangle})d\xi\\
\end{aligned}
$$
We consider two subcases:

{\bf Subcase 1.1: $|y_1'-y_1|\geq 2^{-k}$}

In this situation we have, for any integer $M\geq 0$,
$$
\begin{aligned}
\Big|K_k(y_0,y_1,\ldots, y_m)&-K_k(y_0,y'_1,\ldots, y_m)\Big|\\
&\leq \Big|\int_{\RR^{mn}}a_k(y_0,\xi)\prod_{j= 1}^m e^{i\langle \xi_j, y_0-y_j \rangle}d\xi\Big|\\
&~~+ \Big|\int_{\RR^{mn}}a_k(y_0,\xi)\prod_{j\neq 1}^m e^{i\langle \xi_j, y_0-y_j \rangle}e^{i\langle \xi_1, y_0-y'_1 \rangle}d\xi\Big|\\
&:=E_1+E_2.
\end{aligned}
$$
By integration by part,
 $|y_0-y_1|\approx |y_0-y'_1|\approx \sum_{k,l}|y_k-y_l|$, (\ref{eq1-pseudodiff}) and definition of the class $S^0_{1,0}$ we have
$$
\begin{aligned}
E_1&\leq \sum_{|\alpha|=mn+N} |y_0-y_1|^{-mn-M}\Big|\int_{\RR^{mn}}a_k(y_0,\xi)\partial_{\xi_1}^{\alpha}\prod_{j= 1}^m e^{i\langle \xi_j, y_0-y_j \rangle}d\xi\Big|\\
&\leq \sum_{|\alpha|=mn+N}|y_0-y_1|^{-mn-M}\Big|\int_{\RR^{mn}}\partial_{\xi_1}^{\alpha}a_k(y_0,\xi)\prod_{j= 1}^m e^{i\langle \xi_j, y_0-y_j \rangle}d\xi\Big|\\
&\leq C|y_0-y_1|^{-mn-M} 2^{-k(mn+M-mn)}.
\end{aligned}
$$
Note that by interpolation the inequality above holds for all $M\geq 0$. Taking $M=N+\epsilon$, we have
$$
\begin{aligned}
E_1&\leq C|y_0-y_1|^{-mn-\epsilon-N} 2^{-k(N+\epsilon)}\\
&\leq C\f{|y_1'-y_1|^\epsilon}{(\sum_{k,l}|y_k-y_l|)^{mn+\epsilon}}|y_0-y_1|^{-N}.
\end{aligned}
$$
Hence,
$$
E_1\leq C\f{|y_1'-y_1|^\epsilon}{(\sum_{k,l}|y_k-y_l|)^{mn+\epsilon}}\min\{1,|y_0-y_1|^{-N}\}(2^k|y_1-y_1'|)^{-\epsilon}.
$$
Likewise,
$$
E_2\leq C\f{|y_1'-y_1|^\epsilon}{(\sum_{k,l}|y_k-y_l|)^{mn+\epsilon}}\min\{1,|y_0-y_1|^{-N}\}(2^k|y_1-y_1'|)^{-\epsilon}.
$$
Hence, in this situation, $T_a$ satisfies (H2).

\medskip

{\bf Subcase 1.2: $|y_1'-y_1|< 2^{-k}$}

In this case, by integration by part, we write, for an integer $M\geq 0$,
$$
\begin{aligned}
\Big|K_k(y_0,y_1,\ldots, y_m)&-K_k(y_0,y'_1,\ldots, y_m)\Big| \\
&=\sum_{|\alpha|=M}|y_0-y_1|^{-M}\Big|\int_{\RR^{mn}}a_k(y_0,\xi)(1-e^{i\langle \xi_1, y_1-y_1' \rangle})\partial^\alpha_{\xi_1}\prod_{j= 1}^m e^{i\langle \xi_j, y_0-y_j \rangle}d\xi\\
&=\sum_{|\alpha|=M}|y_0-y_1|^{-M}\Big|\int_{\RR^{mn}}\partial^\alpha_{\xi_1}\Big[a_k(y_0,\xi)(1-e^{i\langle \xi_1, y_1-y_1' \rangle})\Big]\prod_{j= 1}^m e^{i\langle \xi_j, y_0-y_j \rangle}d\xi\\
&=\sum_{|\alpha|+|\beta|=M}|y_0-y_1|^{-M}\Big|\int_{\RR^{mn}}\partial^\alpha_{\xi_1}a_k(y_0,\xi)\partial^\beta_{\xi_1}(1-e^{i\langle \xi_1, y_1-y_1' \rangle})\prod_{j= 1}^m e^{i\langle \xi_j, y_0-y_j \rangle}d\xi.
\end{aligned}
$$
If $|\beta|=0$, $|1-e^{i\langle \xi_1, y_1-y_1' \rangle}|\leq |\xi_1||y_1-y_1'|\leq 2^k|y_1-y_1'|$. Hence
$$
\begin{aligned}
\sum_{|\alpha|=M}|y_0-y_1|^{-M}\Big|\int_{\RR^{mn}}&\partial^\alpha_{\xi_1}a_k(y_0,\xi)(1-e^{i\langle \xi_1, y_1-y_1' \rangle})\prod_{j= 1}^m e^{i\langle \xi_j, y_0-y_j \rangle}d\xi\\
&\leq C|y_0-y_1|^{-M}2^{-k(M-mn-1)}|y_1-y_1'|.
\end{aligned}
$$

If $|\beta|>0$, $|\partial^\beta_{\xi_1} (1-e^{i\langle \xi_1, y_1-y_1' \rangle})|\leq C|y_1-y_1'|^{|\beta|}$. Hence
$$
\begin{aligned}
\sum_{|\alpha|+|\beta|=M, |\beta|>0}|y_0-y_1|^{-M}\Big|&\int_{\RR^{mn}}\partial^\alpha_{\xi_1}a_k(y_0,\xi)\partial^\beta_{\xi_1}(1-e^{i\langle \xi_1, y_1-y_1' \rangle})\prod_{j= 1}^m e^{i\langle \xi_j, y_0-y_j \rangle}d\xi\\
&\leq C\sum_{|\alpha|=M}|y_0-y_1|^{-M}2^{-k(M-|\beta|-mn)}|y_1-y_1'|^{|\beta|}\\
&\leq C|y_0-y_1|^{-M}2^{-k(M-mn-1)}|y_1-y_1'|.
\end{aligned}
$$
Hence, for any integer $M\geq 0$,
$$
\Big|K_k(y_0,y_1,\ldots, y_m)-K_k(y_0,y'_1,\ldots, y_m)\Big|\leq C|y_0-y_1|^{-M}2^{-k(M-mn-1)}|y_1-y_1'|.
$$
By interpolation, this inequality holds for all $M\geq 0$. Taking $M=N+mn+\epsilon$, we have
$$
\begin{aligned}
\Big|K_k(y_0,y_1,\ldots, y_m)&-K_k(y_0,y'_1,\ldots, y_m)\Big|\\
&\leq C|y_0-y_1|^{-N-mn-\epsilon}2^{-k(N+\epsilon-1)}|y_1-y_1'|\\
&\leq C|y_0-y_1|^{-N-mn-\epsilon}|y_1-y_1'|^\epsilon (2^k|y_1-y_1'|)^{1-\epsilon}.
\end{aligned}
$$
Therefore, in this situation, $T_a$ satisfies (H2).

\vskip1cm

{\bf Case 2: $j= 0$}

We can assume that $\max_{k}|y_k-y_1|:=|y_0-y_1|$. We also consider two subcases:

{\bf Subcase 2.1: $|y_0-y'_0|>2^{-k}$}

By the similar argument to that of Subcase 1.1 we get that $T_a$ satisfies (H2).

{\bf Subcase 2.2: $|y_0-y'_0|\leq 2^{-k}$}

In this case, we have
$$
\begin{aligned}
K_k(y_0,y_1,&\ldots, y_m)-K_k(y'_0,y_1,\ldots, y_m)\\
&=\int_{\RR^{mn}}\Big[a_k(y_0,\xi)\prod_{j=1}^m e^{i\langle \xi_j, y_0-y_j\rangle}-a_k(y'_0,\xi)\prod_{j=1}^m e^{i\langle \xi_j, y'_0-y_j\rangle}\Big]d\xi\\
&=\int_{\RR^{mn}}\Big[a_k(y_0,\xi)-a_k(y'_0,\xi)\Big]\prod_{j=1}^m e^{i\langle \xi_j, y_0-y_j\rangle}d\xi\\
&~~~~+\sum_{\ell=1}^m \int_{\RR^{mn}}a_k(y'_0,\xi)\Big[e^{i\langle \xi_\ell, y_0-y_\ell\rangle}-e^{i\langle \xi_\ell, y'_0-y_\ell\rangle}\Big] \prod_{j<\ell}e^{i\langle \xi_j, y'_0-y_j\rangle}\prod_{j>\ell}e^{i\langle \xi_j, y_0-y_j\rangle}  d\xi\\
&=I +\sum_{\ell=1}^mI_\ell.
\end{aligned}
$$
Repeating the arguments as in Subcase 1.2, we get that for any $N\geq 0$
\begin{equation}\label{eq1.1}
\sum_{\ell=1}^m |I_\ell|\leq C|y_0-y_1|^{-N-mn-\epsilon}|y_0-y_0'|^\epsilon (2^k|y_0-y_0'|)^{1-\epsilon}
\end{equation}
for some $\epsilon>0$.

It remains to estimate $I$. We can write, for any integer $M\geq 0$,
$$
\begin{aligned}
|I|&\leq C\sum_{|\alpha|=M}|y_0-y_1|^{-M}\Big|\int_{\RR^{mn}}\Big[a_k(y_0,\xi)-a_k(y'_0,\xi)\Big]\partial_{\xi_1}^\alpha e^{i\langle \xi_1, y_0-y_1\rangle}\prod_{j\neq 1}e^{i\langle \xi_j, y_0-y_j\rangle}d\xi\Big|\\
&\leq C\sum_{|\alpha|=M}|y_0-y_1|^{-M}\Big|\int_{\RR^{mn}}\partial_{\xi_1}^\alpha\Big[a_k(y_0,\xi)-a_k(y'_0,\xi)\Big]\prod_{j= 1}^m e^{i\langle \xi_j, y_0-y_j\rangle}d\xi\Big|
\end{aligned}
$$
where in the last inequality we use integration by part.

By the Mean value Theorem and definition of the class $S^0_{1,\delta}$, we have
$$
\Big|\partial_{\xi_1}^\alpha\Big[a_k(y_0,\xi)-a_k(y'_0,\xi)\Big]\Big|\leq C|y_0-y_0'|2^{-k|\alpha|+k\delta}=C|y_0-y_0'|2^{-kM+k\delta}.
$$
Therefore, for any integer $M\geq 0$,
\begin{equation}\label{eq1.2}
|I|\leq C|y_0-y_1|^{-M}|y_0-y_0'|2^{-k(M-mn)+k\delta}.
\end{equation}
By interpolation again, (\ref{eq1.2}) still holds for any $M\geq 0$.

We now choose $\epsilon>0$ so that $1-\delta>\epsilon$. Taking $M=mn+N+\epsilon$, we have
\begin{equation}\label{eq1.3}
\begin{aligned}
|I|&\leq C|y_0-y_1|^{-(mn+N+\epsilon)}|y_0-y_0'|^\epsilon 2^{-k(N-\delta-\epsilon +1)}(2^k|y-y_0|)^{1-\epsilon}.
\end{aligned}
\end{equation}
The combination of (\ref{eq1.1}) and (\ref{eq1.3}) yields that $T_a$ satisfies (H2).
This completes our proof.
\begin{flushright}
    $\Box$
\end{flushright}
\medskip
From Proposition \ref{pro1-pseudo}, Theorem \ref{thm1} and Theorem \ref{thm2}, we imply the following result.

\begin{thm}\label{thm3}
Let $a\in m$-$S^0_{1,\delta}$, $0\leq \delta<1$. If $T_a$ satisfies (H3), then the following statements hold:

(i) For $1<p_1,\ldots, p_m<\vc$ and $\vec{w}\in A_{\vec{p}}^\vc$, we have
$$
\|T_a(\vec{f})\|_{L^p(\nu_{\vec{w}})}\leq C\prod_{j=1}^m \|f_j\|_{L^{p_j}(w_j)};
$$

(ii) If $1\leq p_1,\ldots, p_m<\vc$ and at least one of the $p_j=1$, then
$$
\|T_a(\vec{f})\|_{L^{p,\vc}(\nu_{\vec{w}})}\leq C\prod_{j=1}^m \|f_j\|_{L^{p_j}(w_j)};
$$

(iii) For any $\vec{b}\in BMO_{\vec{\theta}}$, $1<p_1,\ldots, p_m<\vc$ and $\vec{w}\in A^\vc_{\vec{p}}$, we have
$$
\|(T_a)_{\vec{b}}(\vec{f})\|_{L^p(\nu_{\vec{w}})}\leq C\prod_{j=1}^m \|f_j\|_{L^{p_j}(w_j)}.
$$
\end{thm}

\medskip

We would like to give some relevant comments:

\medskip

(i) It was proved in \cite{GT} that if $T_a$ be a multilinear pseudodifferential operator with the symbol $a$ in the
class $m$-$S^0_{1,\delta}, 0\leq \delta<1$ and all of the transposes $T^{*j}_a$ also have symbols in $m$-$S^0_{1,1}$, then $T_a$ satisfies (H3). Hence, the conclusions in Theorem  \ref{thm3} hold for such a $T_a$.

\medskip

(ii) In particular case when $m=2$, the authors in \cite{BT} proved that if $T_a$ be a multilinear pseudodifferential operator with the symbol $a$ in the
class $2$-$S^0_{1,\delta}, \delta\in [0,1)$ then all of the transposes $T^{*1}_a$ and $T^{*2}_a$ also have symbols in $2$-$S^0_{1,\delta}, \delta\in [0,1)$. Therefore, by the previous remark, the results in Theorem  \ref{thm3} hold for $T_a$.

\medskip

(iii) It can be believed that the obtained results in Theorem  \ref{thm3} still hold for multilinear pseudodifferential operators with the symbols $a\in$ $m$-$S^l_{\rho,\delta}$ with  $0\leq \delta \leq \rho\leq 1,$ $\delta<1, 0<\rho, l<mn(\rho-1)$. These results will be studied in the forth-coming paper.

\end{document}